\documentclass{gtart_h}

\def\ifplaintex{\expandafter\ifx\csname documentclass\endcsname\relax}

\def\gtp{{\mathsurround=0pt\it $\cal G\mskip-2mu$eometry \&\ 
$\cal T\!\!$opology $\cal P\!$ublications}}  

\def\recd{{\small Received:\qua\receiveddate\ifx\reviseddate\relax
\else\qquad Revised:\qua\reviseddate\fi\par}} 

\def\lognumber#1{\def\thelognumber{#1}}
\def\volumenumber#1{\def\thevolumenumber{#1}}
\def\volumeyear#1{\def\thevolumeyear{#1}}
\def\papernumber#1{\def\thepapernumber{#1}}
\def\pagenumbers#1#2{\def\startpage{#1}\def\finishpage{#2}}
\def\published#1{\def\publishdate{#1}}

\def\received#1{\def\receiveddate{#1}}

\def\accepted#1{\def\accepteddate{#1}}

\def\asciiaddress#1{\def\theasciiaddress{#1}}
\def\asciiemail#1{\def\theasciiemail{#1}}

\long\def\asciiabstract#1{\long\def\theasciiabstract{#1}}

\let\\\par\let\thelognumber\relax\let\thevolumenumber\relax
\let\thepapernumber\relax\let\thevolumeyear\relax\let\startpage\relax
\let\finishpage\relax\let\publishdate\relax\let\receiveddate\relax
\let\reviseddate\relax\let\accepteddate\relax\let\theasciititle\relax
\let\theasciiauthors\relax\let\theasciiaddress\relax
\let\theasciiabstract\relax

\let\theasciiemail\relax

\ifplaintex
\font\logobig=cmssbx10 scaled 3836
\font\logomed=cmssbx10 scaled 2557
\else
\font\logobig=cmssbx10 scaled 4200
\font\logomed=cmssbx10 scaled 2800
\fi

\long\def\makeagttitle{   
\count0=\startpage
\agt\hfill      
\hbox to 45truept{\vbox to 0pt{\vglue -13truept{\logomed A\kern -.37em{\logobig 
T}\kern -.38em G}\vss}\hss}
\break
{\small Volume \thevolumenumber\ (\thevolumeyear)
\startpage--\finishpage\nl
Published: \publishdate}

\vglue .25truein

{\parskip=0pt\leftskip 0pt plus
1fil\def\\{\par\smallskip}{\Large\bf\thetitle}\par\medskip} \vglue
0.05truein

{\parskip=0pt\leftskip 0pt plus 1fil\def\\{\par}{\sc\theauthors}
\par\medskip}%
 
\vglue 0.03truein 

{\small\leftskip 25truept\rightskip 25truept{\bf Abstract}\stdspace\theabstract

{\bf AMS Classification}\stdspace\theprimaryclass
\ifx\thesecondaryclass\relax\else; \thesecondaryclass\fi\par
{\bf Keywords}\stdspace \thekeywords\par}\vglue 7truept

}   

\ifplaintex
\font\phead=cmsl9 scaled 950
\font\pnum=cmbx10 scaled 913
\font\pfoot=cmsl9 scaled 950
\headline{\vbox to 0pt{\vskip -4.5mm\line{\small\phead\ifnum
\count0=\startpage ISSN 1472-2739 (on-line) 1472-2747 (printed)
\hfill {\pnum\folio}\else\ifodd\count0\def\\{ }%
\ifx\theshorttitle\relax\thetitle\else\theshorttitle\fi\hfill{\pnum\folio}
\else\def\\{ and }{\pnum\folio}\hfill\ifx\theshortauthors\relax\theauthors
\else\theshortauthors\fi\fi\fi}\vss}}
\footline{\vbox to 0pt{\vglue 0mm\line{\small\pfoot\ifnum\count0=\startpage
\copyright\ \gtp\hfill\else
\agt, Volume \thevolumenumber\ (\thevolumeyear)\hfill\fi}\vss}}
\else
\font=cmsl9 scaled 1050
\font\lnum=cmbx10 
\font=cmsl9 scaled 1050
\makeatletter
\def\@oddhead{{\small\lhead\ifnum\count0=\startpage ISSN 1472-2739 
(on-line) 1472-2747 (printed)\hfill {\lnum\number\count0}\else\ifodd\count0
\def\\{ }\ifx\theshorttitle\relax \thetitle \else\theshorttitle\fi\hfill
{\lnum\number\count0}\else\def\\{ and }{\lnum\number\count0}
\hfill\ifx\theshortauthors\relax 
\theauthors\else\theshortauthors\fi\fi\fi}}\def\@evenhead{\@oddhead}
\def\@oddfoot{\small\lfoot\ifnum\count0=\startpage\copyright\ \gtp\hfill\else
\agt, Volume \thevolumenumber\ (\thevolumeyear)\hfill\fi}
\def\@evenfoot{\@oddfoot}
\makeatother
\fi
\let\maketitlepage\makeagttitle
\let\makeshorttitle\maketitlepage
\let\maketitle\maketitlepage

\newwrite\gtoutfile
\long\gdef\makeheadfile{  
{\def\\{, }\def\s{ }
\immediate\openout\gtoutfile head.xxx
\immediate\write\gtoutfile{Proxy-for: \ifx\theasciiauthors\relax
\theauthors\else\theasciiauthors\fi\s<\ifx\theasciiemail\relax\theemail\else\theasciiemail\fi>}
\immediate\write\gtoutfile{\noexpand\\}
\immediate\write\gtoutfile{Authors: \ifx\theasciiauthors\relax
\theauthors\else\theasciiauthors\fi}
{\def\\{ }\immediate\write\gtoutfile{Title: \ifx\theasciititle\relax
\thetitle\else\theasciititle\fi}}
\immediate\write\gtoutfile{Subj-class: GT or SG, GR etc}
\immediate\write\gtoutfile{MSC-class: \theprimaryclass\ifx\thesecondaryclass\relax\else, \thesecondaryclass\fi}
\immediate\write\gtoutfile{Journal-ref: Algebr. Geom. Topol. \thevolumenumber\s
(\thevolumeyear) \startpage-\finishpage}
\immediate\write\gtoutfile{Comments: Published by Algebraic and
Geometric Topology at}
\immediate\write\gtoutfile{\s\s\s  http://www.maths.warwick.ac.uk/agt/AGTVol\thevolumenumber/agt-\thevolumenumber-\thepapernumber.abs.html}
\immediate\write\gtoutfile{\noexpand\\}
\immediate\write\gtoutfile{}
\ifx\theasciiabstract\relax
\immediate\write\gtoutfile{\theabstract}\else
\immediate\write\gtoutfile{\theasciiabstract}\fi
\immediate\write\gtoutfile{}
\immediate\write\gtoutfile{\noexpand\\}
\immediate\write\gtoutfile{}
\immediate\closeout\gtoutfile}}  

\def\maketitlepage{\makeagttitle\makeheadfile}
\let\makeshorttitle\maketitlepage
\let\maketitle\maketitlepage

\lognumber{67}
\volumenumber{5}
\volumeyear{2005}
\papernumber{67}
\pagenumbers{1655}{1676}
\received{7 November 2005} 
\accepted{10 November 2005}
\published{1 December 2005}

\usepackage{amssymb,epsfig}

\def\fref#1{\hyperlink{#1anchor}{\ref*{#1}}}
\def\figref#1{\hyperlink{#1anchor}{Figure~\ref*{#1}}}
\def\anchor#1{\noindent\hypertarget{#1anchor}{\smash{$\phantom{99}$}}}

\newtheorem{theorem}{Theorem}
\newtheorem{lemma}[theorem]{Lemma}
\newtheorem*{Theorem}{Theorem}
\newtheorem{corollary}[theorem]{Corollary}
\newtheorem{proposition}[theorem]{Proposition}

\theoremstyle{definition}
\newtheorem{example}[theorem]{Example}
\newtheorem{definition}[theorem]{Definition}
\newtheorem{remark}[theorem]{Remark}
\newtheorem{step}{Step}

\newcommand{\bdy}{\mbox{bdy }}
\newcommand{\closure}{\mbox{cl }}
\newcommand{\interior}{\mbox{int }}

\begin{document}

\title[Fundamental groups of subsets of surfaces inject into
their shape groups]{The
fundamental groups of subsets of closed surfaces\\inject into
their first shape groups}
\authors{Hanspeter Fischer\\Andreas Zastrow}

\address{Department of Mathematical Sciences, Ball State 
University\\Muncie, IN 47306, USA}
\secondaddress{Institute of Mathematics, University of
Gda\'nsk\\ul. Wita Stwosza 57, 80-952 Gda\'nsk,  Poland}

\asciiaddress{Department of Mathematical Sciences, Ball State 
University\\Muncie, IN 47306, USA\\and\\Institute 
of Mathematics, University of
Gdansk\\ul. Wita Stwosza 57, 80-952 Gdansk,  Poland}

\asciiemail{fischer@math.bsu.edu, zastrow@math.univ.gda.pl}
\gtemail{\mailto{fischer@math.bsu.edu}{\rm\qua 
and\qua}\mailto{zastrow@math.univ.gda.pl}}

\begin{abstract}
We show that for every subset $X$ of a closed surface $M^2$ and every
$x_0\in X$, the natural homomorphism $\varphi:\pi_1(X,x_0)\rightarrow
\check{\pi}_1(X,x_0)$, from the fundamental group to the first shape
homotopy group, is injective.  In particular, if $X\subsetneq M^2$ is
a proper compact subset, then $\pi_1(X,x_0)$ is isomorphic to a
subgroup of the limit of an inverse sequence of finitely generated
free groups; it is therefore locally free, fully residually free and
residually finite.
\end{abstract}

\asciiabstract{%
We show that for every subset X of a closed surface M^2 and every
basepoint x_0, the natural homomorphism from the fundamental group to
the first shape homotopy group, is injective.  In particular, if X is
a proper compact subset of M^2, then pi_1(X,x_0) is isomorphic to a
subgroup of the limit of an inverse sequence of finitely generated
free groups; it is therefore locally free, fully residually free and
residually finite.}


\maketitle

\section{Introduction}
Understanding the fundamental group of a locally complicated space
can be a difficult endeavor. Indeed, it was shown only recently
 that the fundamental group of the Sierpi\'nski carpet
is not isomorphic to a subgroup of the fundamental group of the
Hawaiian Earrings \cite{eda}. Therefore, properties of such groups
are often deduced indirectly. For example,  the natural
homomorphism $\varphi:\pi_1(X,x_0)\rightarrow
\check{\pi}_1(X,x_0)$, from the fundamental group of any
1-dimensional\linebreak (Hausdorff) compactum $(X,x_0)$ to its
first shape homotopy group, has been shown to be injective
\cite{eda-kawamura}. Consequently, the fundamental group of
any\linebreak 1-dimensional (metric) continuum
 is isomorphic to a subgroup
of the limit of an inverse sequence of finitely generated free
groups; a fact already proved in \cite{curtis-fort} and also in
\cite[Theorem~5.11]{cannon-conner}. Injectivity of the above
homomorphism has also been established for certain fractal-like
trees of manifolds, which need not be semilocally simply connected
at any point \cite{fischer-guilbault}.

 In this article, we
make this same approach available for all subsets of closed
surfaces. Our main result, Theorem~\ref{general}, states that, for
every subset $X$ of a closed surface $M^2$ and every $x_0\in X$,
the naturally induced homomorphism
$\varphi:\pi_1(X,x_0)\rightarrow \check{\pi}_1(X,x_0)$ is
injective. In particular, if $X\subsetneq M^2$ is a proper compact
subset, then $\pi_1(X,x_0)$ is isomorphic to a subgroup of the
limit of an inverse sequence of finitely generated free groups.
Consequently, $\pi_1(X,x_0)$ is locally free, fully residually
free and residually finite (Corollary~\ref{local}).

Another  helpful notion in studying fundamental groups of
complicated spaces is that of being homotopically Hausdorff,
which was introduced in \cite{cannon-conner}.
Proposition~\ref{Hausdorff} below shows that if $X$ is any space
such that  $\varphi:\pi_1(X,x_0)\rightarrow \check{\pi}_1(X,x_0)$
is injective, then $X$ is  homotopically Hausdorff at $x_0$; but
the converse does not hold, in general. Therefore, our main
theorem improves upon \cite[Theorem~3.4]{conner-lamoreaux}, which
states that every subset of the plane is homotopically Hausdorff
at every one of its points.

The structure of this article is as follows:
Section~\ref{definitions} provides a brief review of the
definition of the first shape homotopy group.  In
Section~\ref{body}, we then proceed to prove the main theorem. Our
proofs will use some special properties of subsets of the plane,
which are exported into Section~\ref{tools}. Among them are two
crucial tools from the thesis \cite{zastrowH}:
Theorem~\ref{zastrowA} and Theorem~\ref{zastrowB}. Since these
results are central to our main theorem, we include a complete
proof of Theorem~\ref{zastrowA} and sketch a proof of
Theorem~\ref{zastrowB} in the appendix.

\smallskip

{\bf Acknowledgements}\qua This research was supported by the Faculty
Internal Grants Program of Ball State University as well as  the
University of Gda\'nsk (research grant numbers BW UG
5100-5-0233-2, 5100-5-0096-5). The authors wish to thank these
institutions for their hospitality during their respective visits.
The authors would also like to thank the referee for very helpful
remarks.

\section{The first shape homotopy group} \label{definitions}

\noindent We recall  the  definition of the first shape homotopy
group of a pointed separable metric space $(X,x_0)$. For more
details, see \cite[Chapter~II, \S3.3, Chapter I, \S5.4, and
Appendix~1, \S3]{mardesic-segal}.

Consider the collection $\mathcal C$ of all pairs $({\mathcal U},
\ast)$, where $\mathcal U$ is a locally finite open cover of $X$,
and  $\ast\in \mathcal U$ is a designated member with $x_0\in
\ast$. Then $\mathcal C$ is naturally directed by refinement.
Denote by $(N({\mathcal U}),\ast)$ a geometric realization of the
pointed nerve of $\mathcal U$, that is, a geometric realization of
the abstract simplicial complex $\{\Delta\mid
\emptyset\not=\Delta\subseteq {\mathcal U}, \bigcap_{U\in \Delta}
U \not=\emptyset\}$  with distinguished vertex $\ast$. For every
$({\mathcal U},\ast),({\mathcal V},\ast)\in \mathcal C$ such that
$({\mathcal V},\ast)$ refines $({\mathcal U},\ast)$, choose a
pointed simplicial map $p_{\mathcal U \mathcal V}:(N({\mathcal
V}),\ast)\rightarrow (N(\mathcal U),\ast)$ with the property that
the vertex corresponding to an element $V\in \mathcal V$ gets
mapped to a vertex corresponding to an element $U\in \mathcal U$
with $V\subseteq U$. (Any assignment on the vertices which is
induced by the refinement property will extend linearly.) Then
$p_{\mathcal U \mathcal V}$ is unique up to pointed homotopy and
we denote its pointed homotopy class by $[p_{\mathcal U \mathcal
V}]$. For each $(\mathcal U,\ast)\in \mathcal C$ choose a pointed
map $p_{\mathcal U}:(X,x_0)\rightarrow (N(\mathcal U),\ast)$ such
that $p_{\mathcal U}^{-1}(St(U,N(\mathcal U)))\subseteq U$ for all
$U\in \mathcal U$, where $St(U,N(\mathcal U))$ denotes the open
star of the vertex of $N(\mathcal U)$ which corresponds to $U$.
(For example, define
 $p_{\mathcal U}$ based on a
 partition of unity subordinated to $\mathcal U$.)
Again, such a map $p_{\mathcal U}$ is unique up to pointed
homotopy and we denote its pointed homotopy class by
$[p_{\mathcal U}]$. Then $[p_{\mathcal U \mathcal V}\circ
p_{\mathcal V}]=[p_{\mathcal U}]$. The so-called (pointed)
{\em \v{C}ech expansion}
\[(X,x_0)\stackrel{([p_{\mathcal U}])}{\longrightarrow} (
(N(\mathcal U),\ast), [p_{\mathcal U\mathcal V}],\mathcal C)\] is
an  HPol$_\ast$-expansion, so that we can define the {\em first
shape homotopy group} of~$X$, based at $x_0$, by
\[
\check{\pi}_1(X,x_0)=\lim_{\longleftarrow} \;( \pi_1(N(\mathcal
U),\ast), p_{\mathcal U\mathcal  V\#},\mathcal C).\]

\begin{remark}\label{compactsequence}
In the event that $X$ is compact, other HPol$_\ast$-expansions
are possible. For any inverse sequence
\[(X_1,x_1)\stackrel{f_{1,2}}{\longleftarrow} (X_2,x_2)
\stackrel{f_{2,3}}{\longleftarrow} (X_3,x_3)
\stackrel{f_{3,4}}{\longleftarrow} \cdots\] of pointed compact
polyhedra such that \[(X,x_0)=\lim_{\longleftarrow}
((X_i,x_i),f_{i,i+1}),\] the inverse limit projection
\[(X,x_0)\stackrel{(p_i)}{\longrightarrow}
((X_i,x_i),f_{i,i+1})\] with  $p_i=f_{i,i+1}\circ p_{i+1}$ for
all $i$, is also an HPol$_\ast$-expansion. Therefore,
\[\check{\pi}_1(X,x_0)=\lim_{\longleftarrow}\Big(\pi_1(X_1,x_1)
\stackrel{f_{1,2\#}}{\longleftarrow} \pi_1(X_2,x_2)
\stackrel{f_{2,3\#}}{\longleftarrow} \pi_1(X_3,x_3)
\stackrel{f_{3,4\#}}{\longleftarrow} \cdots\Big).\]
 For example,
an admissible choice is any nested sequence $X_1\supseteq X_2
\supseteq \cdots $ of compact polyhedra, whose intersection is
$X$, where $f_{i,i+1}$ is inclusion and $x_i=x_0$ for all $i$.
If, on the other hand, $X$ is an arbitrary topological space, a
definition of $\check{\pi}_1(X,x_0)$ is based on normal coverings
$(\mathcal U, \ast)$ of $X$, rather than locally finite ones
\cite{mardesic-segal}.
\end{remark}

Since the maps $p_{{\mathcal U}}$ induce homomorphisms
$p_{{\mathcal U}\#}:\pi_1(X,x_0)\rightarrow \pi_1(N({\mathcal
U}),\ast)$ such that  $p_{{\mathcal U}\#}=p_{{\mathcal U \mathcal
V}\#}\circ p_{{\mathcal V}\#}$, whenever $({\mathcal V},\ast)$
refines $({\mathcal U},\ast)$, we obtain an induced homomorphism
\[\varphi:\pi_1(X,x_0)\rightarrow \check{\pi}_1(X,x_0)\] given by
$\varphi([\alpha])=([\alpha_{\mathcal U}])$ where
$\alpha_{\mathcal U}=p_{\mathcal U}\circ \alpha$; or, in the
compact case of Remark~\ref{compactsequence}, by
$\varphi([\alpha])=([\alpha_1],[\alpha_2],[\alpha_3],\dots)$ where
$\alpha_i=p_i\circ\alpha$.

\section{The main theorem}\label{body}
We first establish our result for subsets of the plane:

\begin{theorem}\label{non-comp-planar}
Let $X\subseteq \mathbb{R}^2$ be any  subset  and $x_0\in X$.
Then the natural homomorphism $\varphi:\pi_1(X,x_0)\rightarrow
\check{\pi}_1(X,x_0)$ is injective.
\end{theorem}

Before we present the proof of Theorem~\ref{non-comp-planar}, we
would like to point out that it fails for subsets of
$\mathbb{R}^3$; even for Peano continua in $\mathbb{R}^3$. This is
illustrated in the following example, which  is taken from
\cite[Example~0.12]{zastrow-van-kampen} and based on
\cite[pp.~185--186]{griffiths}.

\begin{example}\label{Hawaiian} Let $H$ be the Hawaiian Earrings.
  That
is, let $H=\bigcup_{k=1}^\infty H_k$ where $H_k=\{(x,y)\in
\mathbb{R}^2\mid x^2+(y-1/k)^2=(1/k)^2\}$. Denote by $C$ the cone
on $H$, with distinguished point $x_0=(0,0)$ in its base~$H$.
Identify two identical copies of $C$ along their distinguished
points  within $\mathbb{R}^3$, and call the resulting space $X$.
Then $\varphi:\pi_1(X,x_0)\rightarrow \check{\pi}_1(X,x_0)$ is
{\em not} injective. Indeed, a loop $\alpha:(S^1,\ast)\rightarrow
(X,x_0)$ which runs through every circle of the two bases exactly
once, always alternating between the two copies of $C$, is not
null-homotopic in $X$, but null-homotopic in every compact
polyhedral approximation of $X$ in $\mathbb{R}^3$.
\end{example}

\begin{remark}
It is well known that the natural homomorphism of
Theorem~\ref{non-comp-planar} is typically not surjective, even
for planar Peano continua. Indeed, such is the case for the
Hawaiian Earrings $(H,x_0)$ of Example~\ref{Hawaiian}.
 To see
why, put $X_i=H_1\cup H_2\cup\cdots \cup H_i$, let
$x_0=x_i=(0,0)$, and define $f_{i,i+1}:X_{i+1}\rightarrow X_i$ by
$f_{i,i+1}(z)=x_0$ for $z\in H_{i+1}$ and $f_{i,i+1}(z)=z$ for
$z\in X_{i+1}\setminus H_{i+1}$. Then $(H,x_0)$ is the limit of
the inverse sequence $((X_i,x_i),f_{i,i+1})$. While this time
$\varphi: \pi_1(H,x_0)\rightarrow \check{\pi}_1(H,x_0)$ is
injective \cite[Corollary~1.2]{eda-kawamura}, it is not
surjective:
 let $l_i:(S^1,\ast)\rightarrow (H_i,x_0)$ be a fixed
homeomorphism and, following an idea of Griffiths
\cite[p.~185]{griffiths-two}, consider for each $i$ the element
\[g_i=[l_1][l_1][l_1]^{-1}[l_1]^{-1}
[l_1][l_2][l_1]^{-1}[l_2]^{-1}
[l_1][l_3][l_1]^{-1}[l_3]^{-1}\cdots
[l_1][l_i][l_1]^{-1}[l_i]^{-1}\] of $\pi_1(X_i,x_i)$. Then the
sequence $(g_i)$ is an element of the group
$\check{\pi}_1(H,x_0)$ which is clearly not in the image of
$\varphi:\pi_1(H,x_0)\rightarrow \check{\pi}_1(H,x_0)$.
\end{remark}

\begin{proof}[Proof of Theorem~\ref{non-comp-planar}]
 Let $\alpha:(S^1,\ast)\rightarrow (X,x_0)$ be a map with the
 property that
$[\alpha]\not=1\in\pi_1(X,x_0)$. We wish to show that
$\varphi([\alpha])\not=1\in \check{\pi}_1(X,x_0)$. To this end,
define \[Y=\alpha(S^1)\cup \bigcup\{U\mid U \mbox{ is a (path)
component of } \mathbb{R}^2\setminus \alpha(S^1) \mbox{ with }
U\subseteq X\}.\] We now apply the results of Section~\ref{tools}.
By Lemma~\ref{Peano}, $Y$ is a Peano continuum. Since $x_0\in
Y\subseteq X$, we have $[\alpha]\not=1\in\pi_1(Y,x_0)$. Let $L$
and $(L_i)_i$ be as in
 Theorem~\ref{zastrowA}, applied to the planar Peano continuum $Y$ with designated
 point~$x_0$, and arrange for $z_i\not\in X$ for all $i$. Since the
inclusion $(Y,x_0)\hookrightarrow (L,x_0)$ is a pointed homotopy
equivalence, we have
 $[\alpha]\not=1\in\pi_1(L,x_0)$.
 By Theorem~\ref{zastrowB}, there
 is an $n$ such that
 $[\alpha]\not=1\in\pi_1(L_n,x_0)$.
Let $\{D_1,D_2,\dots,D_{n-1}\}$ be the closures of the bounded
components of $\mathbb{R}^2\setminus L_n$. Then  $z_i\in
(\mbox{int } D_i)\setminus X$ for all $1\leqslant i\leqslant n-1$.
Since
 $[\alpha]\not=1\in\pi_1(L_n,x_0)$,
 Lemma~\ref{inner-region} implies that $\alpha$ is not null-homotopic
in  $\mathbb{R}^2\setminus \{z_1, z_2, \dots, z_{n-1}\}$. Choose a
finite open cover $\mathcal V$ of $\mathbb{R}^2\setminus \{z_1,
z_2, \dots, z_{n-1}\}$, such that $p_{{\mathcal
V}\#}:\pi_1(\mathbb{R}^2\setminus \{z_1, z_2, \dots,
z_{n-1}\},x_0)\rightarrow \pi_1(N({\mathcal V}),\ast)$ is an
isomorphism. Define $\mathcal U=\{V\cap X\mid V\in {\mathcal
V}\}$. Then  ${\mathcal U}$ is a finite open cover of $X$. Notice,
that $N({\mathcal U})$ is a subcomplex of $N({\mathcal V})$, and
that the two maps $p_{\mathcal V}|_X:X \rightarrow N({\mathcal
V})$ and $p_{\mathcal U}:X\rightarrow N({\mathcal U})\subseteq
N({\mathcal V})$ are  homotopic as maps into $N({\mathcal V})$,
because they are contiguous. (That is, for every $x\in X$,
$p_{\mathcal V}(x)$ and $p_{\mathcal U}(x)$ lie in a common
simplex of $N({\mathcal V})$.) If $\varphi([\alpha])=1\in
\check{\pi}_1(X,x_0)$, then  $p_{{\mathcal U}\#}([\alpha]) =1 \in
\pi_1(N({\mathcal U}),\ast)$. This would imply that
 $p_{{\mathcal V}\#}([\alpha]) =1
\in  \pi_1(N({\mathcal V}),\ast)$, which is not the case.
\end{proof}

\begin{theorem}\label{general}
Let  $M^2$ be a closed surface and $X\subseteq M^2$ any subset.
Then for every $x_0\in X$, the natural homomorphism
$\varphi:\pi_1(X,x_0)\rightarrow \check{\pi}_1(X,x_0)$ is
injective.
\end{theorem}

\begin{proof}
We will argue by way of contradiction and assume, contrary to the
assertion, that there is a continuous function
$\alpha:(S^1,\ast)\rightarrow (X,x_0)$ with the property that
$[\alpha]\not=1\in \pi_1(X,x_0)$ and
$\varphi([\alpha])=1\in\check{\pi}_1(X,x_0)$. Consider the
following commutative diagram of naturally induced homomorphisms:
\begin{eqnarray*}
\pi_1(X,x_0)& \stackrel{incl_{\#}}{\longrightarrow} & \pi_1(M^2,x_0) \\
\downarrow \;\;\;\;\;&\circlearrowleft & \;\;\;\;\downarrow \\
\check{\pi}_1(X,x_0)& \stackrel{incl_{\#}}{\longrightarrow} &
\check{\pi}_1(M^2,x_0) \end{eqnarray*} Since the righthand
vertical homomorphism is an isomorphism, we conclude that
$[\alpha]=1\in \pi_1(M^2,x_0)$ and $X\not=M^2$. Replacing $M^2$,
if necessary, we may therefore assume that $M^2$ has universal
covering $q:(\mathbb{R}^2,\tilde{x}_0)\rightarrow (M^2,x_0)$ so
that $\alpha$ lifts to a map $\tilde{\alpha}:(S^1,\ast)\rightarrow
(\mathbb{R}^2, \tilde{x}_0)$ with $q\circ \tilde{\alpha}=\alpha$.
Let us define $\bar{X}=q^{-1}(X)$. Then $\tilde{x}_0\in
\bar{X}\subseteq \mathbb{R}^2$. Consider
 the natural homomorphism $\psi:\pi_1(\bar{X},\tilde{x}_0)\rightarrow
\check{\pi}_1(\bar{X},\tilde{x}_0)$.

We {\em claim} that
$\psi([\tilde{\alpha}])=1\in\check{\pi}_1(\bar{X},\tilde{x}_0)$.
To verify this, let $\mathcal U$ be any cover of $\bar{X}$ by open
sets of $\bar{X}$. Choose a sufficiently fine cover $\mathcal V$
of $\bar{X}$ by open sets of $\bar{X}$ such that (i) $\mathcal V$
refines $\mathcal U$; (ii) ${\mathcal C}=\{q(V)\mid V\in {\mathcal
V}\}$ is a cover of $X$ by open sets of $X$;  and (iii) there is a
covering map $q'$  (when restricted to path components of
$N({\mathcal V})$) in a commutative diagram:
\begin{eqnarray*}
\bar{X}& \stackrel{p_{\mathcal V}}{\longrightarrow} & N({\mathcal V}) \\
 ^{\tilde{\alpha}}\!\!\! \nearrow \;\;\;  \downarrow^q\!
  &\circlearrowleft & \;\;\;\downarrow^{q'}\\ S^1 \stackrel{\alpha}{\longrightarrow}  X&
\stackrel{p_{\mathcal C}}{\longrightarrow} & N({\mathcal C})
\end{eqnarray*}
 Since
$\varphi([\alpha])=1\in\check{\pi}_1(X,x_0)$, we have
$p_{{\mathcal C}\#}([\alpha])=1\in\pi_1(N({\mathcal C}),\ast)$.
Since $p_{\mathcal V}\circ \tilde{\alpha}$ is a lift of
$p_{\mathcal C}\circ \alpha$, we have $p_{{\mathcal
V}\#}[\tilde{\alpha}]=1\in \pi_1(N({\mathcal V}),\ast)$. Hence,
every open cover $\mathcal U$ of $\bar{X}$ is refined by an open
cover $\mathcal V$ of $\bar{X}$ such that $p_{{\mathcal
V}\#}[\tilde{\alpha}]=1\in \pi_1(N({\mathcal V}),\ast)$.
Consequently,
$\psi([\tilde{\alpha}])=1\in\check{\pi}_1(\bar{X},\tilde{x}_0)$,
as claimed.

However, $[\tilde{\alpha}]\not=1\in\pi_1(\bar{X},\tilde{x}_0)$,
because  $[\alpha]\not=1\in \pi_1(X,x_0)$. This contradicts
Theorem~\ref{non-comp-planar}. \end{proof}

Recall that a compact metric space of trivial shape is called {\em
cell-like}. Recall also that a compact metric space is cell-like
if and only if it contracts within every ANR that it is embedded
in. Every cell-like subset $X$ of a 2-dimensional manifold $M^2$
is {\em cellular} in $M^2$, that is, $X$ has arbitrarily small
neighborhoods in $M^2$, which are homeomorphic to open disks.
While the path components of a cell-like planar continuum need
not be contractible  (e.g. for the Knaster continuum), we record
as a corollary to Theorem~\ref{general} the following
``folklore'' fact:

\begin{corollary}\label{cell-like}
If $X$ is a cell-like subset of a 2-dimensional manifold $M^2$ and
$x_0\in X$, then $\pi_1(X,x_0)=1$.
\end{corollary}

\begin{proof}
Since the first shape homotopy group of a cell-like space is
trivial, this follows directly from Theorem~\ref{general}.
\end{proof}

 Recall that a group $G$ is said to be {\em
locally free} if each of its finitely generated subgroups is free.
A group $G$ is called {\em fully residually free} if for every
finite subset $\{g_1, g_2,\dots ,g_n\}\subseteq G\setminus\{1\}$
there is a homomorphism $\phi:G\rightarrow F$ to a free group $F$
such that $\phi(g_i)\not= 1$ for all $i\in\{1,2,\dots,n\}$.
Finally, $G$ is said to be {\em residually finite} if for every
$g\in G\setminus\{1\}$ there is a homomorphism $\phi:G\rightarrow
H$ to a finite group $H$ such that $\phi(g)\not=1$.

 Some of the algebraic implications of
Theorem~\ref{general} are collected in the following:

\begin{corollary}\label{local}
Let $X$ be a proper compact subset of a  2-dimensional manifold
$M^2$  and
 $x_0\in X$. Then $\pi_1(X,x_0)$ is isomorphic to a subgroup of the limit of an inverse
 sequence of finitely generated free groups. In particular,  $\pi_1(X,x_0)$ is
locally free, fully  residually free  and residually
 finite.
\end{corollary}

\begin{proof}
Notice that we can express $(X,x_0)$ as the intersection of a
nested sequence $(X_i,x_0)$ of compact 2-dimensional PL manifold
neighborhoods of $X$ in $M^2$, with $\partial X_i\not=\emptyset$
for all $i$. Since each of the groups $\pi_1(X_i,x_0)$ is finitely
generated and free, it follows from Theorem~\ref{general} and
Remark~\ref{compactsequence} that $\pi_1(X,x_0)$ is isomorphic to
a subgroup of the limit of an inverse
 sequence of finitely generated free groups.

A proof of the fact that every finitely generated subgroup of the
 limit of an inverse sequence of finitely generated free groups is free can be found in
 \cite[Theorem~1]{curtis-fort}.
 Using the projection maps of this inverse limit to its free terms, we see that
 $\pi_1(X,x_0)$ is  fully residually free. The first proof of the
 fact that free groups are residually finite was given by Schreier
 \cite[p.~97]{chandler-magnus}. (This proof can also be found in
 \cite[p.~11]{cohen}.) It follows that $\pi_1(X,x_0)$ is
 residually finite, since it is fully residually free.
\end{proof}

Following \cite[Definition~5.2]{cannon-conner}, we call a space
$X$ {\em homotopically Hausdorff at the point} $x_0\in X$, if for
every $g\in \pi_1(X,x_0)\setminus\{1\}$ there is an open set
$U\subseteq X$ with $x_0\in U$ such that there is no loop
$\alpha:(S^1,\ast)\rightarrow (U,x_0)$ with $[\alpha]=g$. If $X$
is homotopically Hausdorff at every one of its points, then $X$
is said to be {\em homotopically Hausdorff}.  In
\cite[Theorem~3.4]{conner-lamoreaux} it is shown that all planar
sets are homotopically Hausdorff. Theorem~\ref{general} above is a
strengthening of this statement, because we have the following:

\begin{proposition}\label{Hausdorff}
Let $X$ be a topological space and  $x_0\in X$.  If the natural
homomorphism  $\varphi:\pi_1(X,x_0)\rightarrow
\check{\pi}_1(X,x_0)$ is injective, then $X$ is homotopically
Hausdorff at $x_0$. However, the converse need not hold.
\end{proposition}

\begin{proof}
An element $g\in \pi_1(X,x_0)$, witnessing a violation of the
definition of homotopically Hausdorff at $x_0$,  will be mapped
to the trivial element of the fundamental group of the pointed
nerve of every open cover of $X$, contradicting the injectivity
of $\varphi:\pi_1(X,x_0)\rightarrow \check{\pi}_1(X,x_0)$.
Example~\ref{sine} below exhibits a space which is homotopically
Hausdorff and for which $\varphi$ is not injective.
\end{proof}

\begin{example} \label{sine} This space was previously studied in
\cite{fischer-guilbault-workshop}. Let
\[Y=\{(x,y,z)\in \mathbb{R}^3\mid z=0, 0<x\leqslant 1, y=\sin (1/x)\}\cup
(\{0\}\times[-1,1]\times\{0\})\] be the ``topologist's sine
curve''. Define $Y_i=Y\cup ([0,1/i]\times [-1,1]\times\{0\})$.
Let $X$ and $X_i$ be the subsets of $\mathbb{R}^3$ obtained by
revolving $Y$ and $Y_i$ about the $y$-axis, respectively, and let
$f_{i,i+1}:X_{i+1}\hookrightarrow X_i$ be inclusion. Then $X$ is
the limit of the inverse sequence $(X_i,f_{i,i+1})$. If we take
$x_0=(1,\sin 1,0)$, then $\pi_1(X,x_0)$ is infinite cyclic, while
$\check{\pi}_1(X,x_0)$ is trivial. We observe that
$\varphi:\pi_1(X,x_0)\rightarrow \check{\pi}_1(X,x_0)$ is not
injective, while $X$ is homotopically Hausdorff.
\end{example}

Combining Theorem~\ref{general} with Proposition~\ref{Hausdorff}
also yields a proof of the following:

\begin{corollary}
Every subset $X$ of a closed surface is homotopically Hausdorff.
\end{corollary}

\section{Some special properties of subsets of the plane}\label{tools}

 \begin{definition}\label{innerdef}
Consider a map $\alpha:S^1\rightarrow \mathbb{R}^2$ and let $O$ be
the unbounded (path) component of $\mathbb{R}^2\setminus
\alpha(S^1)$. We will call $\mathbb{R}^2\setminus O$ the {\em
inner region} of $\alpha$ and denote it by $I(\alpha)$.
\end{definition}

The inner region $I(\alpha)$ of any map $\alpha:S^1\rightarrow
\mathbb{R}^2$
 is a Peano continuum. This follows from the following lemma.

\begin{lemma}\label{Peano}
Let $X\subseteq Y$ be two Peano continua, i.e., two compact,
connected and locally (path) connected metric spaces.
 Let $\mathcal C$ be
the collection of (path) components of $Y\setminus X$ and
$\mathcal C'\subseteq \mathcal C$. Then $Z=X\cup \bigcup \mathcal
C'$ is also a Peano continuum.
\end{lemma}

\begin{proof} By assumption, $(Y,d)$ is a metric space.
We will verify that $Z$ is (i)~compact, (ii) connected, and (iii)
locally path connected.

(i)\qua Since each element of $\mathcal C$ is open in $Y$, then
$Z=Y\setminus \bigcup [\mathcal C\setminus \mathcal C']$ is a
closed subset of the compact space $Y$. Hence $Z$ is compact.

(ii)\qua Since $X$ is path connected, it will suffice to show that for
every $z\in U\in \mathcal C'$ there is an $x\in X$ and a path
$p:[0,t]\rightarrow Z$ from $p(0)=z$ to $p(t)=x$. \linebreak To
this end, pick any $y\in X$ and choose a path $q:[0,1]\rightarrow
Y$ from $q(0)=z$ to $q(1)=y$. Let $t=\min [q^{-1}(Y\setminus U)]$.
Then $p=q|_{[0,t]}$ and $x=p(t)$ will work.

(iii)\qua Let $x\in Z$. Note that every element of $\mathcal C'$ is
an open subset of the locally path connected space  $Y$.  Hence,
in order to show that $Z$ is locally path connected at~$x$, we may
assume that $x\in X$. Let $\epsilon>0$ be given. We will find a
$\delta>0$ such that every point of $Z\cap N(x,\delta)$ can be
joined to $x$ by a path that lies in $Z\cap N(x,\epsilon)$. Here,
we adopt the usual notation $N(y,\eta)=\{z\in Y \mid
d(z,y)<\eta\}$. Since $X$ is locally path connected, we can
choose an $\eta\in(0,\epsilon)$ such that every point of $X\cap
N(x,\eta)$ can be joined to $x$ by a path that lies in $X\cap
N(x,\epsilon)$. Since $Y$ is also locally path connected, we can
choose a $\delta\in(0,\eta)$ such that every point of
$N(x,\delta)$ can be joined to $x$ by a path that lies in
$N(x,\eta)$.

Now, let $y\in Z\cap N(x,\delta)$. We may assume that $y\not\in
X$. Choose $U\in \mathcal C'$ with $y\in U$. Choose a path
$p:[0,1]\rightarrow N(x,\eta)$ from $p(0)=y$ to $p(1)=x$. Let
$t=\min[p^{-1}(Y\setminus U)]$. Then  $p|_{[0,t]}$ is a path from
$y$ to $p(t)\in X\cap N(x,\eta)$ that lies in $U\cap N(x,\eta)$,
and therefore lies in $Z\cap N(x,\epsilon)$. Since there is
another path in $X\cap N(x,\epsilon)$ from $p(t)$ to $x$, we are
done.
\end{proof}

\begin{lemma}\label{inner-region} Let any set $X\subseteq \mathbb{R}^2$ and any map
$\alpha:S^1\rightarrow X$ be given. Let $I(\alpha)$ be as in
Definition~\ref{innerdef}.
 If $\alpha:S^1\rightarrow X$
is null-homotopic, then so is $\alpha:S^1\rightarrow X\cap
I(\alpha)$.
\end{lemma}

\begin{proof}
Let us identify $\mathbb{R}^2$ with $\mathbb{S}^2\setminus
\{\infty\}$ by means of stereographic projection and write
 $I(\alpha)=\mathbb{S}^2\setminus U$, where $U$ is the (path)
component of $\infty$ in $\mathbb{S}^2\setminus \alpha(S^1)$.
Observe that $U$ is a simply connected open subset of
$\mathbb{S}^2$. (For if there were a homotopically non-trivial
loop in $U$, which we might as well assume to be piecewise linear,
it would have to contain a simple closed curve which is
non-trivial in $U$. However, such a non-trivial simple closed
curve of $U$ would have to separate points of $I(\alpha)$. This
would contradict the path connectedness of $I(\alpha)$.)
 Moreover, we may assume that the complement of
$U$ in $\mathbb{S}^2$, namely $I(\alpha)$, contains more than one
point.
 Hence, by the Riemann Mapping Theorem,
there is a biholomorphic map $g:\mathbb{D}\rightarrow U$ whose
domain is the open unit disk $\mathbb{D}\subseteq \mathbb{C}$.
Since $I(\alpha)$ is locally path connected, there is a map
$\bar{g}:\bar{\mathbb{D}}\rightarrow \mathbb{S}^2$, whose domain
is the closed unit disk $\bar{\mathbb{D}}$, such that
$\bar{g}|_{\mathbb{D}}=g$ and $\bar{g}(\partial \bar{\mathbb{D}})=
\bdy U$ \cite[Theorem~2.1]{pommerenke}. Put $d=g^{-1}(\infty)\in
\mathbb{D}$ and choose a continuous retraction
$r:\bar{\mathbb{D}}\setminus \{d\}\rightarrow
\partial \bar{\mathbb{D}}$ with $r|_{\partial
\bar{\mathbb{D}}}=id$. Since
$\bar{g}|_{\bar{\mathbb{D}}\setminus\{d\}}:\bar{\mathbb{D}}\setminus
\{d\}\rightarrow [\closure U]\setminus \{\infty\}$ is a quotient map,
$r$  induces a continuous retraction  $R:[\closure
U]\setminus\{\infty\}\rightarrow \bdy U\subseteq I(\alpha)$ such
that $R|_{\mbox{{\small bdy $U$}}}=id$,  which we can extend  to
a continuous retraction $R:\mathbb{R}^2\rightarrow I(\alpha)$ by
defining $R(x)=x$ for all $x\in I(\alpha)$.

Now, if $H:S^1\times[0,1]\rightarrow X$ is a homotopy contracting
$\alpha$ to a point within $X$, that is, $H(s,0)=\alpha(s)$ for
all $s\in S^1$ and $H(S^1,1)$ is a singleton, then the homotopy
$R\circ H:S^1\times[0,1]\rightarrow I(\alpha)$ contracts $\alpha$
to a point within $X\cap I(\alpha)$. Indeed, if $H(s,t)\in
I(\alpha)$, then $R(H(s,t))=H(s,t)\in X$, and if $H(s,t)\in U$,
then we have $R(H(s,t))\in \bdy U\subseteq \alpha(S^1)\subseteq
X$, also.
\end{proof}

\begin{definition}
A sequence $(C_i)$ of subsets of a metric space is called a {\em
null-sequence} if for every $\epsilon>0$,
$\mbox{diam}(C_i)\leqslant \epsilon$ for all but finitely many
$i$.
\end{definition}

The following theorem is based on \cite{zastrowH}.

\begin{theorem}\label{zastrowA}
Let $Y\subseteq \mathbb{R}^2$ be a Peano continuum. Fix any
$x_0\in Y$. List the bounded (path) components of
$\mathbb{R}^2\setminus Y$ in a (possibly finite) pairwise disjoint
sequence $C_1, C_2, C_3, \dots $. From each $C_i$ let one point
$z_i$ be given.
 Then
there is a sequence $L_1, L_2, L_3,\dots$ of subsets of
$\mathbb{R}^2$ such that
 $L_1$ is a closed disk and each $L_{i+1}$ is obtained from
$L_i$ by removing the interior of a closed disk $D_i$ in the
interior of $L_i$ with the following properties: $z_i\in \mbox{int
} D_i$ for all $i$;\linebreak $(D_i)$ forms a null-sequence;
 $Y\subseteq L_i$ for all $i$ and, setting $L=\bigcap_{\;i} L_i$,
 the inclusion $(Y,x_0) \hookrightarrow (L,x_0)$ is a pointed homotopy
 equivalence.
\end{theorem}

The proof of Theorem~\ref{zastrowA} begins with the following lemma.

\begin{lemma}\label{null}
Let $Y\subseteq \mathbb{R}^2$ be a Peano continuum. Then any
pairwise distinct sequence of bounded (path) components of
$\mathbb{R}^2\setminus Y$ is a null-sequence.
\end{lemma}

\begin{proof}
Suppose, to the contrary, that there is an $\epsilon>0$ and
infinitely many distinct bounded components $C_1, C_2, C_3, \dots$
of $\mathbb{R}^2\setminus Y$  with $\mbox{diam}( C_i)> \epsilon$
for all $i$.\linebreak For each $i\in\mathbb{N}$ choose points
$x_i=(a_i,b_i)\in C_i$ and $y_i=(c_i,d_i) \in C_i$ with
$d(x_i,y_i)>\epsilon$. Since $\bigcup_i C_i$ is bounded, so are
the sequences $(x_i)$ and $(y_i)$. Passing to subsequences, we may
therefore assume that there are $x,y\in\mathbb{R}^2$ such that
$(x_i)$ converges to $x$ and $(y_i)$ converges to $y$. Then
$d(x,y)\geqslant \epsilon$. Changing our coordinate system, if
necessary, we may assume that $x=(-2,0)$ and $y=(2,0)$.\linebreak
For $j\in\{-1,0,1\}$, consider the three vertical lines
$L_j=\{(j,r)\mid r\in \mathbb{R}\}$. Discarding finitely many
indices, we may assume that $a_i<-1$ and $c_i>1$ for all
$i\in\mathbb{N}$. Since every $C_i$ is a path connected open
subset of the plane containing points to the left of $L_{-1}$ as
well as points to the right of $L_1$, we may choose for each
$i\in\mathbb{N}$ a polygonal arc $\alpha_i$ in $C_i$ such that
$\alpha_i\cap L_{-1}=\{s_i\}$, $\alpha_i\cap L_1=\{t_i\}$,
$\partial \alpha_i=\{s_i,t_i\}$, and $\alpha_i$ is in general
position with respect to $L_0$. Passing to another subsequence, we
may assume that the sequence $(s_i)$ is monotone, say (strictly)
increasing. Then the sequence $(t_i)$ also increases (strictly),
because the arcs $\alpha_1, \alpha_2, \alpha_3,\dots$ are pairwise
disjoint. For each fixed $i\in\mathbb{N}$, consider the closed
disk $D_i$, bounded by the simple closed curve consisting of
$\alpha_i$, $\alpha_{i+1}$, and the two vertical closed intervals
$[s_i,s_{i+1}]$ and $[t_i,t_{i+1}]$. Since $L_0$ separates
$[s_i,s_{i+1}]$ from $[t_i,t_{i+1}]$, there is a vertical interval
$[u_i,v_i]\subseteq L_0\cap D_i$ such that $u_i$ lies in one of
the two arcs $\alpha_i$ or $\alpha_{i+1}$ and $v_i$ lies in the
other arc. Now, $[u_i,v_i]$ is an arc connecting the two distinct
complementary components $C_i$ and $C_{i+1}$ of $Y$ in
$\mathbb{R}^2$. Hence, we can choose a point $w_i\in (u_i,v_i)\cap
Y$. Since $(w_i)$ is a bounded sequence in $L_0$, we may pass to a
subsequence and assume that $(w_i)$ converges to some $w\in L_0$.
In fact, since $Y$ is closed and each $w_i\in Y$, we must have
$w\in Y$. By construction, $D_i\cap D_{i+1}=\alpha_{i+1}$ and
$D_i\cap D_j=\emptyset$ for all $j>i+1$. Also,
 the points $w_i$ lie in the pairwise disjoint interiors of the
 disks~$D_i$. Consequently, $w\not\in D_i$ for all $i\in\mathbb{N}$. We
conclude that every path in $Y$ connecting $w$ to $w_i$ must meet
either $[s_i,s_{i+1}]\subseteq L_{-1}$ or $[t_i,t_{i+1}]\subseteq
L_1$, because $\alpha_i\cap Y=\emptyset$ for all $i\in
\mathbb{N}$. Consequently, $Y$ is not locally path connected at
$w$, which is a contradiction.
\end{proof}

\begin{proof}[Proof of Theorem~\ref{zastrowA}]
As in the proof of Lemma~\ref{inner-region}, we identify
$\mathbb{R}^2$ with $\mathbb{S}^2\setminus\{\infty\}$. We put
$z_0=\infty$ and let $C_0$ denote the (path) component of $z_0$ in
$\mathbb{S}^2\setminus Y$. Note that $\mathbb{S}^2\setminus
C_0=Y\cup \bigcup_{i\geqslant 1} C_i$ is a Peano continuum,  by
Lemma~\ref{Peano}. Therefore, as in the proof of
Lemma~\ref{inner-region}, we obtain a biholomorphic map
$g:(\mathbb{D},\ast)\rightarrow (C_0,z_0)$ with continuous
extension $\bar{g}:\bar{\mathbb{D}}\rightarrow \mathbb{S}^2$ such
that \linebreak $\bar{g}(\partial \mathbb{D})=\bdy C_0$. Choose a
closed disk $D$ with $\ast\in \interior D\subseteq D \subseteq
\mathbb{D}$ and put $D_0=g(D)$. Then $z_0\in \interior
D_0\subseteq D_0\subseteq C_0$. Choose a  deformation retraction
$H$ of $\bar{\mathbb{D}}\setminus \interior D$ onto $\partial
\mathbb{D}$, that is, choose a map
$H:\left(\bar{\mathbb{D}}\setminus \interior
D\right)\times[0,1]\rightarrow \bar{\mathbb{D}}\setminus \interior
D$ with $H(x,0)=x$ for all $x\in \bar{\mathbb{D}}\setminus
\interior D$,  $H(x,1)\in \partial \mathbb{D}$ for all $x\in
\bar{\mathbb{D}}\setminus \interior D$, and $H(x,t)=x$ for all
$x\in\partial \mathbb{D}$ and all $t\in [0,1]$. Since
\[\bar{g}|_{\bar{\mathbb{D}}\setminus \mbox{\scriptsize int}\,
D}\times id: \left(\bar{\mathbb{D}}\setminus \interior
D\right)\times[0,1]\rightarrow \left(\closure C_0 \setminus
\interior D_0\right)\times [0,1]\] is a quotient map, the
deformation retraction $H$ induces a deformation retraction $F_0$
of  $\closure C_0\setminus \interior D_0$ onto $\bdy C_0$.

 If we choose closed disks $D_i$ with
$z_i\in\interior D_i\subseteq D_i\subseteq C_i$ for $i\geqslant 1$
in the same way, we obtain deformation retractions $F_i$ of
$\closure C_i\setminus \interior D_i$ onto $\bdy C_i$ for all
$i\geqslant 0$.

Define $L_1=\mathbb{S}^2\setminus \interior D_0$,
$L_{i+1}=L_i\setminus \interior D_i$ and $L=\bigcap_{\;i} L_i$.
 Since $F_i(x,t)=x$ whenever $x\in \bdy C_i$ and since $C_0, C_1, C_2,
 \dots$ are pairwise disjoint, we may
 combine these deformation retractions
into one  function $F:L\times [0,1]\rightarrow L$ by defining
$F(x,t)=F_i(x,t)$ if $x\in \closure C_i$ for some $i\geqslant 0$,
and $F(x,t)=x$ otherwise. Then
 $F(x,0)=x$ for all $x\in L$, $F(x,1)\in Y$ for all $x\in L$, and $F(x,t)=x$
for all $x\in Y$ and all $t\in [0,1]$. In order to complete the
proof of the theorem, we only need to establish  continuity of
$F$.

Let $(x,t)\in L\times[0,1]$ be given. Since each $F_i$ is
continuous, we may assume that $x$ is not contained in any of the
open sets $C_i$. Then $F(x,t)=x$. Let $\epsilon>0$ be given. By
Lemma~\ref{null}, we may choose $N\in \mathbb{N}$ such that
$\mbox{diam}(\closure C_i)=\mbox{diam}(C_i)\leqslant \epsilon/2$
for all $i>N$. Define $I=\{i\in\{0,1,2,\dots ,N\}\mid x\in
\closure C_i\}$. Choose $\delta>0$ such that (i)
$d(F(x,t),F(y,s))=d(F_i(x,t),F_i(y,s))<\epsilon$ whenever
$(y,s)\in L\times[0,1]$ with $y\in\closure C_i$ for some $i\in I$,
$d(x,y)<\delta$ and $|s-t|<\delta$; (ii) $y\not\in \closure C_i$
whenever $i\in\{0,1,2,\dots,N\}\setminus  I$ and $d(x,y)<\delta$;
and (iii) $\delta<\epsilon/2$.
 Let $(y,s)\in L\times [0,1]$ with
$d(x,y)<\delta$ and $|s-t|<\delta$ be given. If $y\not\in\closure
C_i$ for any $i\geqslant 0$, then $F(y,s)=y$ so that
$d(F(x,t),F(y,s))=d(x,y)<\delta<\epsilon/2<\epsilon$, by (iii). We
may therefore assume that $y\in \closure C_i$ for some $i\geqslant
0$. If $i>N$, then\linebreak
$d(F(x,t),F(y,s))=d(x,F(y,s))\leqslant
d(x,y)+d(y,F(y,s))<\delta+\epsilon/2<\epsilon/2+\epsilon/2$,
because both $y$ and $F(y,s)$ are contained in $\closure C_i$ and
$\mbox{diam}(\closure C_i)\leqslant \epsilon/2$. Hence, we may
assume that $0\leqslant i\leqslant N$. Since $d(x,y)<\delta$ and
$y\in\closure C_i$, then $i\in I$ by (ii). Consequently,
$d(F(x,t),F(y,s))<\epsilon$ by (i).
\end{proof}

The following theorem can be extracted from \cite{zastrowH}, where
its development is a sideline of a comprehensive study of higher
dimensional homotopy groups of planar sets. However, it is not
explicitly stated in \cite{zastrowH} and  its proof is scattered
throughout the lengthy document. We will therefore sketch a proof
in the appendix.

\begin{theorem}\label{zastrowB}
 Let $L_1, L_2, L_3,\dots$ be a sequence of subsets of $\mathbb{R}^2$
such that
 $L_1$ is a closed disk and each $L_{i+1}$ is obtained from
$L_i$ by removing the interior of a closed disk $D_i$ in the
interior of $L_i$ such that $(D_i)$ forms a null-sequence.
 If $x_0\in L=\bigcap_i L_i$, then the natural homomorphism
$\varphi:\pi_1(L,x_0)\rightarrow \check{\pi}_1(L,x_0)$ is
injective.
\end{theorem}

\begin{remark}
While it is known \cite[Corollary~1.2]{eda-kawamura} that
$\varphi:\pi_1(X,x_0)\rightarrow \check{\pi}(X,x_0)$ is injective
for 1-dimensional continua $X$, a space $L$ as described in
Theorem~\ref{zastrowB} need not be homotopy equivalent to a
1-dimensional metric space (see \cite[\S 5]{cannon-conner-zastrow}
and \cite{karimov}).
\end{remark}

\section*{Appendix}
\addcontentsline{toc}{section}{Appendix}

The purpose of this appendix is to sketch a proof of
Theorem~\ref{zastrowB}.

\setcounter{step}{-1}

\begin{step} {\sl Reduction of Theorem~\ref{zastrowB} to a theorem about Sierpi\'nski-like
spaces.}

For $i\in \mathbb{N}$, consider the collection $\mathcal
Q_i^{\ast\ast}$ of all subsquares of $[0,1]\times[0,1]$ of the
form $[\frac{2k-1}{3^i},\frac{2k}{3^i}]\times
[\frac{2m-1}{3^i},\frac{2m}{3^i}]$ with $k,m\in\mathbb{N}$. Let
$\mathcal Q^\ast_i$ be the collection of all elements of
${\mathcal Q}^{\ast\ast}_i$ which are not contained in any element
of  ${\mathcal Q}^{\ast\ast}_s$ with  $s<i$. \figref{modified}
(left) depicts the union of the collections ${\mathcal Q}^\ast_1,
{\mathcal Q}^\ast_2$ and ${\mathcal Q}^\ast_3$ in
$[0,1]\times[0,1]$. Let $S$ be any space obtained from
$[0,1]\times [0,1]$ by removing all, some, or none of the
interiors of the squares in $\bigcup_{\; i}{\mathcal Q}^\ast_i$.
Specifically, for each $i\in\mathbb{N}$, choose any subset
${\mathcal Q_i}\subseteq {\mathcal Q^\ast_i}$. Inductively, define
$S_0=[0,1]\times[0,1]$ and let $S_{i+1}$ be obtained from $S_i$ by
removing the interiors of the elements of ${\mathcal Q}_i$, which
we shall informally refer to as {\em holes}. Then  $S=\bigcap_{\;
i} S_i$ will be called a {\em Sierpi\'nski-like} space with
defining sequence $(S_i)$.
\begin{figure}[ht!]\anchor{modified}
 \begin{center}
\epsfig{file=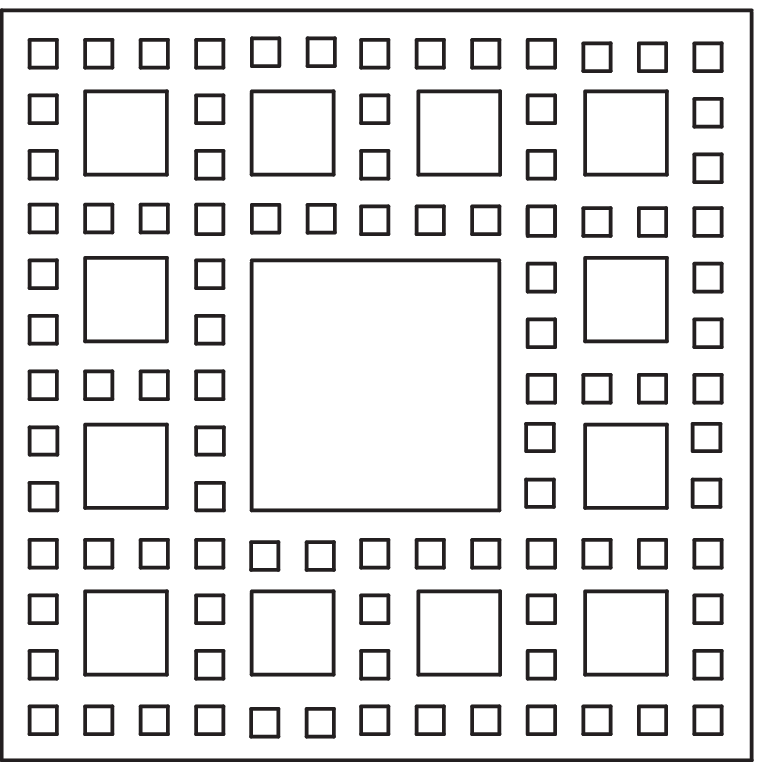,height=2in,width=2in} \hspace{0.5in}
\epsfig{file=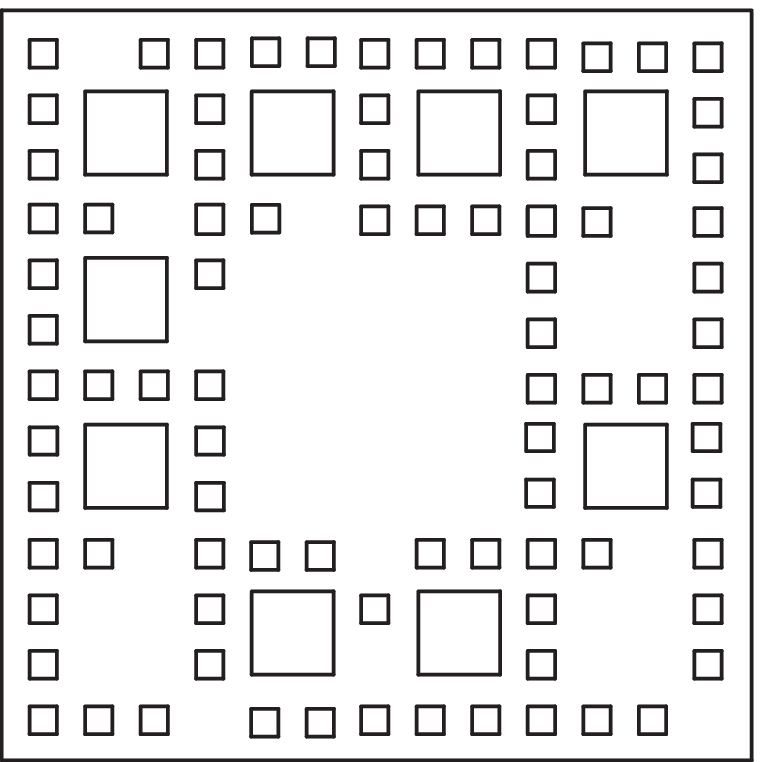,height=2in,width=2in}\hspace{0.5in}
\end{center}
 \caption{A Sierpi\'nski-like pattern;
 complete or incomplete}\label{modified}\end{figure}

Note that if $\mathcal Q_i=\mathcal Q_i^\ast$ for all $i$,  the
resulting space $S$, which we will denote by $S^\ast$, is
homeomorphic to the Sierpi\'nski carpet. This follows from the
topological characterization of the Sierpi\'nski carpet given in
\cite[Theorems~3 and 4]{whyburn}: it is the unique planar
1-dimensional Peano continuum without local cut points. It also
follows that if $L$, $(L_i)$ and $(D_i)$ are as in
Theorem~\ref{zastrowB}, then $L$ is homeomorphic to some
Sierpi\'nski-like space~$S$. Indeed, we may extend the sequence
$(D_i)$ to a dense null-sequence of pairwise disjoint disks
$(D^\ast_i)$ in the interior of $L_1$. It follows from elementary
considerations that removing the interiors of all $D^\ast_i$ from
$L_1$ leaves us with a planar 1-dimensional Peano continuum
$L^\ast$ without local cutpoints, which consequently is
homeomorphic to $S^\ast$.\linebreak Observe that the boundary
curves of the disks whose interiors are removed in either version
of the Sierpi\'nski carpet, $L^\ast$ or $S^\ast$, are
characterized as the unique non-separating simple closed curves.
Therefore, since $L$ is obtained by removing the interiors of only
some of the disks $D^\ast_i$ from $L_1$, namely those of the disks
$D_i$, the homeomorphism between $L^\ast$ and $S^\ast$ can be
extended to a homeomorphism of $L$ onto~$S$.

It therefore suffices to prove the following result.

\begin{Theorem} Let $S$ be a Sierpi\'nski-like  space with defining sequence
$(S_i)$ and $\alpha:S^1\rightarrow S$  a loop which is
null-homotopic in $S_i$ for every $i\in\mathbb{N}$. Then $\alpha$
is null-homotopic in $S$.\end{Theorem}
\end{step}

We will call the closure of a component of $S_i\cap
\left([0,1]\times (\frac{2m-1}{3^i},\frac{2m}{3^i})\right)$ a {\em
horizontal corridor} of $S_i$ (on the $m$-th {\em stratum}) and
denote the collection of all horizontal corridors of $S_i$ by
 ${\mathcal N^h_i}$. Analogously,  we will call the closure of a
component of $S_i\cap
\left((\frac{2k-1}{3^i},\frac{2k}{3^i})\times[0,1]\right)$ a {\em
vertical corridor} of $S_i$ (on the $k$-th {\em stratum}) and
denote the collection of all vertical corridors of $S_i$  by
 ${\mathcal N^v_i}$. If  $[a,b]\times [c,d]\in {\mathcal N}^h_i$ is a horizontal
 corridor, we will refer to $[a,b]\times (c,d)$ as its {\em inner
 region}, to $[a,b]\times\{c\}$ and $[a,b]\times\{d\}$ as its
 {\em boundary lines}, and to $[a,b]\times\{\frac{c+d}{2}\}$ as its {\em center line}.
  Similar terminology will be applied to
 vertical corridors, but in the other direction.

\begin{step}
{\sl The idea is to encode how $\alpha$ traverses the network of
corridors ${\mathcal N}_i={\mathcal N}^h_i\cup {\mathcal N}^v_i$
of every approximation level $S_i$, by associating to $\alpha$ a
sequence of (cyclic) words $\omega=(\omega_1, \omega_2,
\omega_3,\dots)$, each of which is composed of letters that record
which corridor was traversed in what direction, unless $\alpha$
does not cross any corridor of ${\mathcal N}_i$, in which case we
record the empty word for $\omega_i$ on that level.}

 \begin{figure}[ht!]\anchor{word1}
 \parbox{5in}{ \vspace{-.5in} \hspace{.1in}
\epsfig{file=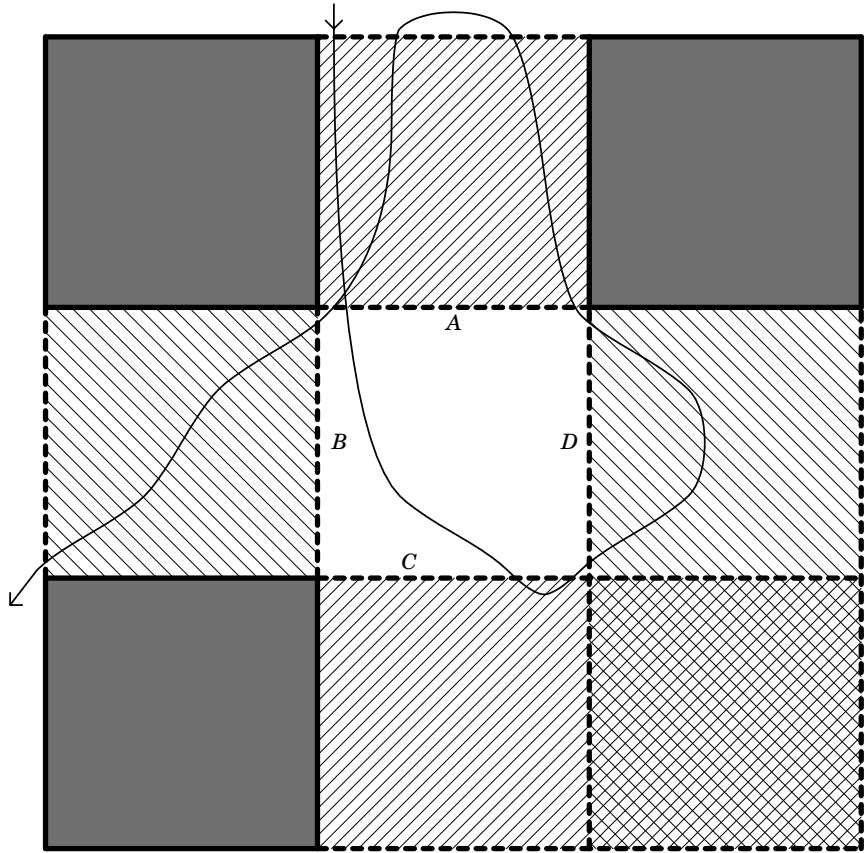,height=5in,width=5in}}

 \caption{$\omega_2=A^{-1}A\;A^{-1}B^{-1}\cdots$}\label{word1}\end{figure}
 Figures~\fref{word1} and \fref{word2} show a
$[\frac{1}{9},\frac{4}{9}]\times[\frac{1}{9},\frac{4}{9}]$-detail
of how such a {\em word sequence} might be formed through two
levels of approximation: $\omega_2=A^{-1}A\;A^{-1}B^{-1}\cdots $
and $\omega_3=a^{-1}_1a^{-1}_3j^{-1}h\;
d_3\;l\;d^{-1}_1a_4\;a_2\;e^{-1}a^{-1}_1a^{-1}_3b^{-1}_3j^{-1}b^{-1}_2\cdots$.

Before we discuss this coding procedure  in more detail in
Step~\ref{code}, we would like to point to the heart of the
problem at hand: consider the loop shown in \figref{cancel},
which runs through a detail of \figref{word2}. This loop is
null-homotopic in $S_3$ and \figref{cancel} depicts a possible
cancellation diagram. Now, the hole on the left hand side of
\figref{cancel} is an element of ${\mathcal Q}_3$. Therefore,
if we consider the same loop in $S_2$ instead of $S_3$,
drastically different cancellations appear plausible, including
those that might not be feasible in $S_3$.
  However, at level $i=2$, we do not know whether or not
such a hole might appear in $S_3$, or if instead a smaller version
of it might appear in some $S_i$ with $i>3$, or even on which side
of the figure it might appear.
 Therefore, a specific null-homotopy of $\alpha$ in
$S_i$, which might have seemed like a natural choice at the time,
will in general be useless for finding a null-homotopy of $\alpha$
in $S$, as it might not even yield a null-homotopy of $\alpha$ in
$S_j$ if $j>i$. The word-calculus of \cite{zastrowH} offers a
solution to this complication.
\begin{figure}[ht!]\anchor{word2}
 \vspace{-.4in}

 \parbox{5in}{\hspace{.1in}
\epsfig{file=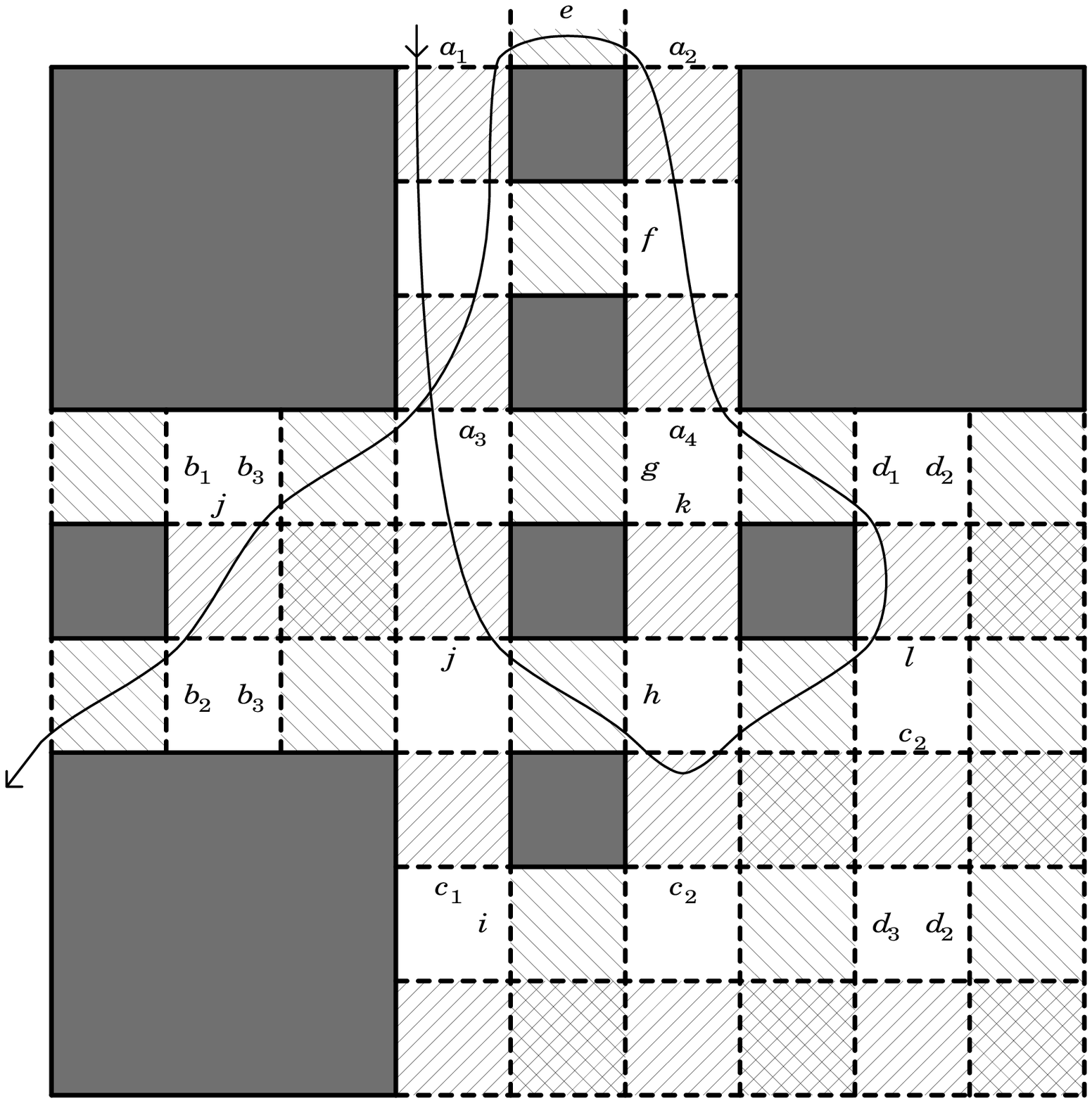,height=5in,width=5in}}

 \caption{$\omega_3=a^{-1}_1a^{-1}_3j^{-1}h\;d_3\;l\;d^{-1}_1a_4\;a_2\;
 e^{-1}a^{-1}_1a^{-1}_3b^{-1}_3j^{-1}b^{-1}_2\cdots$}\label{word2}\end{figure}
\end{step}

\begin{step} \label{code}  {\sl We now describe a  procedure for obtaining a suitable
word sequence.}

 For each $i\in
\mathbb{N}$, consider the collection ${\mathcal I}^h_i$ of all
closed subintervals $[a,b]$ of $S^1$ for which
$\alpha([a,b])\subseteq C$ for some $C\in {\mathcal N}^h_i$, with
$\alpha(a)$ on one boundary line of $C$, $\alpha(b)$ on the other
boundary line of $C$, and such that $\alpha((a,b))$ lies in the
interior region of $C$.

 Then ${\mathcal I}^h_i$ consists of finitely many intervals
with pairwise disjoint interiors. Moreover, we observe the
following relationships between ${\mathcal I}^h_i$ and the
corresponding collections ${\mathcal I}^h_{i-1}$ and ${\mathcal
I}^h_{i+1}$ on its neighboring levels:
\begin{itemize}\item[(R1)] Consider any fixed horizontal stratum
$E$ of level $(i-1)$ and its two substrata $E_1$ and $E_2$ of
level $i$.  If $[a,b]$ and $[c,d]$ are any intervals  of
${\mathcal I}^h_i$ with $\alpha([a,b])\subseteq E_1$ and
$\alpha([c,d])\subseteq E_2$, such that $[a,b]$ and $[c,d]$ are
adjacent in ${\mathcal I}^h_i$ with this property,
 then $[a,d]\in {\mathcal
I}^h_{i-1}$; \item[(R2)] Conversely, for every $[a,d]\in {\mathcal
I}^h_i$, there are $b,c\in [a,d]$ such that both $[a,b],[c,d]\in
{\mathcal I}^h_{i+1}$.
\end{itemize}

Considering the vertical corridors instead, we arrive at a similar
collection~${\mathcal I}^v_i$.

 Form a (cyclic) ordered
tuple $([a_1,b_1], [a_2,b_2], [a_3,b_3],\dots)$ by  merging the
elements of ${\mathcal I}^h_i$ and ${\mathcal I}^v_i$, while
maintaining the relative natural order on $S^1$ within each
collection; though allowing any relative order between two
intersecting intervals $[a, b]\in{\mathcal I}^h_i$ and $[a',
b']\in {\mathcal I}^v_i$ which are associated to crossing
corridors. The word $\omega_i$ is then formed by replacing each
interval $[a,b]$ of this sequence by the label of the corridor
into which it maps under $\alpha$, along with an exponent
indicating positive or negative direction.

Note that the word sequence $(\omega_i)$ contains all relevant
information about the loop~$\alpha$. Indeed, if we were to
recreate $\alpha$ from the combinatorial data contained in the
word sequence $(\omega_i)$ by forming  polygonal loops $\alpha_i$
in $S_i$, which are in general position with respect to the
boundary lines of $\mathcal N_i$ and which spell out the word
$\omega_i$ by canonically traversing the prescribed corridors,
then the sequence $(\alpha_i)$ would  converge uniformly to a
reparametrization of $\alpha$.
\end{step}\vspace{-5pt}

\begin{step}
{\sl Since $\alpha:S^1\rightarrow S_i$
is null-homotopic for every $i\in\mathbb{N}$, we can transform
each word $\omega_i$ to the empty word by symbolic cancellations
and transpositions, if we allow letters to commute whenever they
correspond to crossing corridors.}

Indeed, we can homotope $\alpha:S^1\rightarrow S_i$ to a polygonal
loop $\alpha_i':S^1\rightarrow S_i$ which is in general position
with respect to the center lines of all corridors of ${\mathcal
N}_i$, while maintaining the sequence $([a_1,b_1], [a_2,b_2],
[a_3,b_3],\dots)$ for the word $\omega_i$, and crossing the center
lines only as often as needed to spell the word $\omega_i$. We
then may contract $\alpha_i':S^1\rightarrow S_i$ with a homotopy
$H'_i:D^2\rightarrow S_i$ (also in general position) such that
$H'_i|_{S^1}=\alpha'_i$. Consider the
 preimage of the center
lines of all corridors of ${\mathcal N}_i$ in $D^2$, under the
homotopy $H'_i$. It is the union of finitely many arcs, whose
endpoints lie in $S^1$, and finitely many circles. Intersections
among
 these arcs and circles
  can occur only if one of them maps into
 a horizontal center line and the other one into a vertical center line.
 We can eliminate the circles, one innermost circle at a time, by
moving the homotopy across the appropriate center line. (Note that
this can be done without disturbing the pairwise disjoint arc
components of the preimage of the center lines of the other
direction, which might span across these circles.) The endpoints
of the remaining arcs lie in distinct elements of $([a_1,b_1],
[a_2,b_2], [a_3,b_3],\dots)$, which correspond to letters of the
same corridor but have opposite exponents. (To see why the
exponents cannot agree, lift the entire null-homotopy to the
universal cover of the at most twice punctured unit square;
punctured by the hole(s) at the end of that corridor. An innermost
mismatching arc would, when following the lift of $\alpha_i'$
between its endpoints, contradict the fact that each component of
the preimage of the center line separates this cover.)
 We therefore arrive at a pattern of arcs in $D^2$, which pairs up
cancelling letters of $\omega_i$ across commuting letters as
claimed.

In fact,
 we can
shorten the word in every step of this transformation, so long as
we allow for cancellations to take place across commuting
variables, without formally performing any transpositions.
\end{step}

\begin{step}\label{extend}
{\sl If we repeat the above argument with the combined center
lines of two adjacent levels ${\mathcal N}_i$ and ${\mathcal
N}_{i-1}$, we see that every cancellation facility of $\omega_i$
thus obtained induces at least one cancellation facility of
$\omega_{i-1}$, via the relationships $(R1)$ and $(R2)$ above.
Moreover, we can select a reduction scheme for the entire word
sequence $(\omega_i)$, which is coherent between levels.}

To illustrate the idea, consider the cancellation diagram depicted
in \figref{cancel} for the word $\omega_3=d_3\, d_2\,
l\,d_2^{-1}l^{-1}d_3^{-1}k\,d_1\,
l^{-1}d_2\,l\,d_2^{-1}d_1^{-1}k^{-1}$. It induces cancellations
for the corresponding word $\omega_2=DD^{-1}DD^{-1}$ of the form
$(DD^{-1})(DD^{-1})$, rather than the competing cancellations
$D(D^{-1}D)D^{-1}$. In order to obtain a coherent reduction
scheme for the entire word sequence $(\omega_i)$, we need to find
cancellation diagrams for the individual words $\omega_i$, which
induce each other, from one level to the previous one, in the
above sense.
\begin{figure}[ht!]\anchor{cancel}
\center
\parbox{5in}{\hspace{.2in}
\epsfig{file=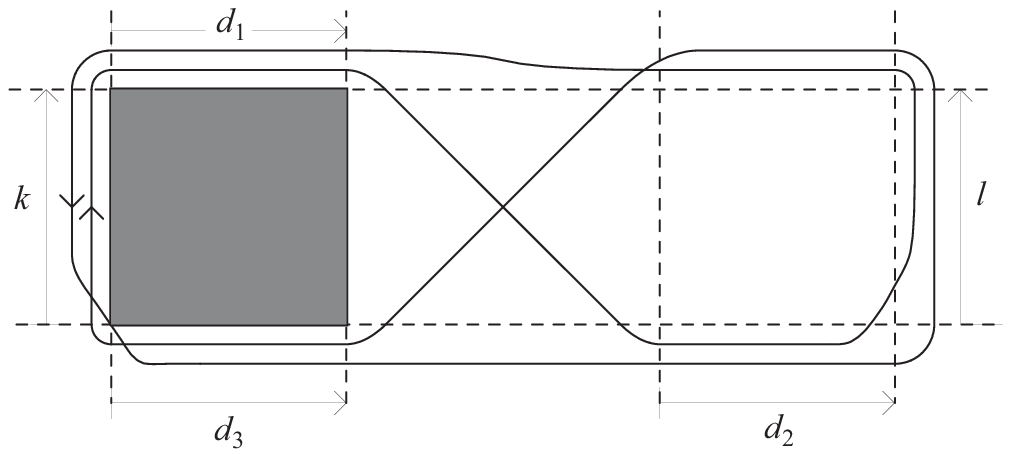}}\vspace{-1.3in}

\parbox{5in}{\center \epsfig{file=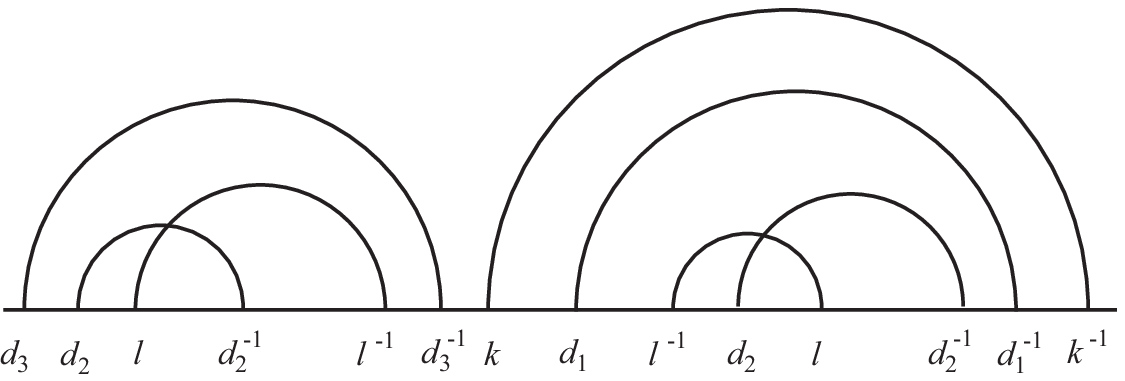}}
\caption{A transformation of the word $\omega_3$ to the empty
word}\label{cancel}\end{figure}

Having reduced this problem to a purely combinatorial issue, we
can find a coherent scheme without much difficulty by considering
the following underlying graph structure. We regard each specific
cancellation of each word $\omega_i$ as one vertex of a graph
whose (directed) edges indicate that a certain cancellation of
$\omega_i$ induces a certain cancellation of $\omega_{i-1}$. One
can now proceed as in K\"onig's  classical method \cite[VI.\S2
Satz~6]{koenig}
 for
 finding an infinite ray in a finitely
branching infinite  rooted tree: beginning at the root-vertex,
select one of the finitely many adjacent vertices behind which
there are infinitely many vertices and construct the infinite ray
 iteratively by repeating this step.
Put differently,
  the collection of cancellation patterns for
the word sequence $(\omega_i)$ can be organized into an inverse
system of nonempty finite sets. The inverse limit of this system
is nonempty and any element of this limit provides a coherent set
of cancellation patterns. (See \cite[Prop. 2.9 -- Cor.
2.11]{morgan-morrison}, where the latter argument is used in a
similar situation.)
\end{step}

\begin{step}\label{cellulation} {\sl The remainder of the proof consists of recording the
resulting cancellation diagram for each level $i\in\mathbb{N}$ and
 turning it into a homotopy
$H_i:D^2\rightarrow S_i$, such that the sequence $(H_i)$ converges
uniformly to a homotopy $H:D^2\rightarrow S$ with
$H|_{S^1}=\alpha$.}

 The transformation information of
\figref{cancel}, for example, can be transcribed into the disk
$D^2$, whose boundary $S^1$ is the domain of $\alpha$, featuring
the intervals ${\mathcal I}^h_i$ and ${\mathcal I}^v_i$. The
result is shown in \figref{cells}, where for each line of the
cancellation diagram of \figref{cancel} we placed a band
between the corresponding intervals of $S^1$. (Recall that an
interval of ${\mathcal I}^h_i$ may overlap with an interval
of~${\mathcal I}^v_i$.)\linebreak

 We thus obtain a sufficiently coherent
sequence
 of cellulations of $D^2$ by bands. Because our transformations
 consist of
 cancellations,
  each such cellulation can
 be represented by (hyperbolic) straight line segments
 whose boundaries lie in~$S^1$.

\begin{figure}[ht!]\anchor{cells}
\center
\parbox{2in}{\epsfig{file=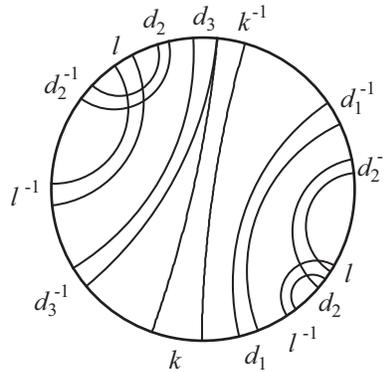,height=2in,width=2in}}
\caption{A cellulation of the disk induced by Figure~\ref{cancel}
}\label{cells}\end{figure}

For a fixed $i\in\mathbb{N}$, we  define $H_i:D^2\rightarrow S_i$ as follows.
 Begin by declaring that  $H_i|_{S^1}=\alpha$. On the two
boundary lines of a band which is associated to a horizontal
letter, we define the $y$-coordinates to be constant. Similarly,
on the two boundary lines of a band which is associated to a
vertical letter, we define the $x$-coordinates to be constant.
Consequently, $H_i$ is defined on the zero-skeleton of our
cellulation. The missing coordinate for a point on the
one-skeleton is now obtained by a linear interpolation of the
corresponding coordinates of the endpoints of this one-cell. By
construction, when $H_i$ is restricted to the boundary of a
two-cell, we obtain a null-homotopic path in $S_i$, yielding the
missing values for $H_i$ in a canonical way. In order to describe
 how, we will refer to a square
$[\frac{t}{3^i},\frac{t+1}{3^i}]\times
[\frac{s}{3^i},\frac{s+1}{3^i}]$ of $S_i$ as being of type $0, 1$,
or $2$, according to how many corridors of ${\mathcal N}_i$
contain its interior. Since $H_i$ maps the boundary lines of each
band into the corresponding corridor of ${\mathcal N}_i$, the
crossing region of two bands, for example, can be mapped into one
type-$2$ square of $S_i$, which is the intersection of the
corresponding corridors of ${\mathcal N}_i$. Similar arguments
apply to the remaining types of two-cells: they can either be
mapped into (i) regions formed by one square of type 1 along with
 at
most two adjacent squares of type 2; or (ii) regions formed by one
center square of type~0 along with at most eight adjacent squares
of type 1 or 2. (Observe that the extended regions described in
(i) and (ii) might be needed to accommodate the images of those
two-cells, which intersect $\partial D^2$. Indeed, as shown in
Figures~\fref{word1} and~\fref{word2}, our loop might enter adjacent
corridors without crossing them, thus not recording the
corresponding letters in the word of this level and hence not
contributing corresponding bands to the cellulation.)
\end{step}

\begin{step}
{\sl We claim that $||H_i(z)-H_{i+1}(z)||\leqslant (6/3^i)$ for
every $i\in \mathbb{N}$ and for every $z\in D^2$. This will imply
that the sequence $(H_i:D^2\rightarrow S_i)$ converges uniformly,
and since $S=\bigcap_{\; i}S_i$, its limit will be a
 homotopy $H:D^2\rightarrow S$ with $H|_{S^1}=\alpha$, thus
completing the proof.}

Let any $i\in \mathbb{N}$ and any $z\in D^2$ be given. Denoting
$H_i(z)=(x_i,y_i)$ and $H_{i+1}(z)=(x_{i+1},y_{i+1})$, it suffices
to argue that $|x_i-x_{i+1}|<(3/3^i)$ and \linebreak
$|y_i-y_{i+1}|<(3/3^i)$. We will only focus on the
$y$-coordinates, since the argument for the $x$-coordinates is
entirely analogous.

For the moment, consider only cellulation bands of level $i$, and
only those which correspond to horizontal letters. Such bands do
not intersect. Moreover, the (constant) $y$-values of $H_i$ on the
boundary lines of these bands provide an estimate for the
$y$-values of $H_i$ on the regions in between: by construction,
$H_i$ never maps such regions across more than three adjacent
horizontal sectors of the form
$[0,1]\times[\frac{s}{3^i},\frac{s+1}{3^i}]$ (see (i) and (ii) of
Step~\ref{cellulation}), so that the $y$-values of $H_i$ on one
such region never differ by more than $3/3^i$.

The same (even stronger) estimates hold for the $y$-values of
$H_{i+1}$ in regions between boundary lines of cellulation bands
of level $(i+1)$, which correspond to horizontal letters.

We claim that such estimates also hold across the two levels $i$
and $(i+1)$. For this to work, it is essential to have selected
coherent cellulations between levels, because two
  bands from  adjacent levels, both corresponding to horizontal
  letters,
 do not randomly cross. In fact, only one single, very benign,
 intersection phenomenon between their boundary lines is possible:

Consider  a loop running through \figref{word2}, whose word
 $\omega_3=a_3a_1a_1^{-1}a_1a_1^{-1}a_3^{-1}$
cancels as
 $a_3(a_1a_1^{-1})(a_1a_1^{-1})a_3^{-1}$ and whose word
 $\omega_2=AA^{-1}$ cancels as $(AA^{-1})$. Let $R$
 be the band representing the cancellation of $A$ with $A^{-1}$ and let
 $r_1$, $r'_1$ and $r_3$ be the bands representing the cancellations of
 $a_1$ with $a_1^{-1}$,  $a_1$ with
 $a_1^{-1}$, and  $a_3$ with
 $a_3^{-1}$, respectively. (See \figref{turning}.) As described in (R1) and (R2) of
 Step~\ref{code}, the bands $r_1$ and
 $r_3$ originate on either side of one end of the band $R$.
 While $r_3$, as one would expect, lies completely inside of $R$, sharing one
 boundary line with $R$, $r_1$ (as well as $r'_1$) turns and
 leaves $R$ by crossing that boundary line of $R$ at which it
 originated.

\begin{figure}[ht!]\anchor{turning}
\begin{center}
\epsfig{file=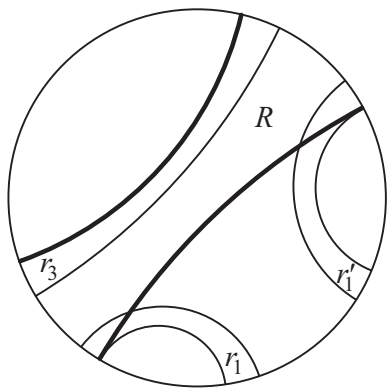,height=2in,width=2in}
\end{center}\caption{Cellulation bands of two adjacent levels --- a worst case scenario}\label{turning}
 \end{figure}

 This is indeed the only type of
 intersection of such bands, because their corresponding arcs
  in the preimage of the homotopy of Step~4, which correlates the cancellation
 diagrams of level $i$ and $(i+1)$, do not
 intersect at all.

 Since there are
 no further intersection phenomena,
  a straightforward region by region analysis of the common refinement of the
 two cellulations on level~$i$ and level~$(i+1)$, by bands corresponding to
 horizontal letters only, therefore  yields
 $|y_i-y_{i+1}|<3/3^i$. \qed

\end{step}

\Addresses\recd

\end{document}